\documentclass[12pt]{article}

\usepackage[T2A]{fontenc}
\usepackage[utf8]{inputenc}

\usepackage{amsfonts}
\usepackage{amssymb}
\usepackage{amsmath}
\usepackage{amsthm}

\def \le {\leqslant}
\def \ge {\geqslant}

\theoremstyle{plain}
 
\begin{document}
\begin{Large}
 \centerline{On an example by Poincar\'{e} and sums with Kronecker sequence
 }
\end{Large}
\vskip+0.5cm
\begin{large}
\centerline{Nikolay Moshchevitin
%\footnote{This research has been financed by the Russian Science Foundation grant 22-21-00079, https://rscf.ru/project/22-21-00079/} 
 }

\end{large}
 
\vskip+1cm
\begin{small}
 {\bf Abstract:}
 This paper is motivated by recent papers by L. Colzani and A. Kochergin. First of all we discuss an example by Poincaré related to sums of the type 
$ \sum_{k=0}^{t-1} f(\alpha k +x)$
where $f$ is a continuous function and $\alpha$ is an irrational number. Then we deal with its generalisations related to multidimensional Kronecker sequence.

 \vskip+0.3cm
 {\bf AMS Subject Classification:} 11K31, 42B99.
 
  \vskip+0.3cm
 {\bf Keywords:} Kronecker sequence, Diophantine Approximation, Diophantine exponents, quasi-periodic functions.
 \end{small}
 \vskip+1cm

 This short and simple communication is motivated by recent papers  \cite{A,Co,K} as well as by the discussion on Poincar\'{e}'s and Bohl's results in \cite{Kz}, Chapter 8.
 We give a brief analysis of an example constructed  by Poincar\'{e}
 and its recent generalisations. Most of the  constructions under consideration are well-known.  In this note, we just wanted to bring all the results  together  and  give a general and improved 
 multi-dimensional formulation of a recent result by A. Kochergin \cite{K}, prove non-existence of a universal continuous function and discuss some of the related results in terms of Diophantine Approximation.  In particular, in our opinion smoothness results involving Diophantine exponents $\omega$, $\hat{\omega}$  and
 $\hat{\lambda}$ had never been documented before.
 
    \vskip+0.3cm
    
 The structure of the paper is as follows. In  Section 1 we briefly describe the original example due to Poincar\'{e}. Section 2 is devoted to a brief survey of  related results for  Birkhoff sums with Kronecker sequence. In Sections 3 and 4 we recall recent results by Kochergin and Colzani correspondingly. Section 5 deals with all necessary objects  related to Diophantine Approximation.
 The main results of the present paper are formulated and proved  in Section 6 (Theorem 1), Section 7 (Theorems 2 and 3) and Section 8 (Theorems 4 and 5).
 
   \vskip+0.3cm
   
    {\bf 1. An example by Poincar\'{e}.}

   \vskip+0.3cm

 Let $A>0$.
 In 1885 Poincar\'{e}  \cite{P1} considered  the following examples of trigonometric series
$$
F_1(t) = \sum_{k=0}^\infty A^{k}\sin \frac{t}{2^k}\,\,\,\,
\text{and}
\,\,\,\,\,
 F_2(t) = 
 \sum_{k=0}^\infty (-A)^{k}\sin \frac{t}{2^k}.
$$
 These series  converge absolutely for $ A< 2$, but for $ A > 1$ they do not converges uniformly in $ t\in \mathbb{R}$. Poincar\'{e}  showed in the case of non-uniform convergence  these 
 series as functions in $t$ may be unbounded.
 Namely, for any $ t_0 \in (0,\sqrt{3}/2)$ put $ B = t_0 -\frac{4}{3}t_0^3>0 $. Suppose that 
 $1+\frac{1}{2B+1} < A< 2$ and define $ h = \frac{2B}{2-A}-\frac{1}{A-1}>0$.
 Poincar\'{e}  proved that  
 for the series $F_1(t)$, because of inequalities
 $$
 F_1(2t) > AF_1(t) -1,\,\,\,\, \forall \, t
 $$
 and
 $$
 F(t) > \frac{B}{1-\frac{A}{2}}\,\,\,\,\,\text{for}\,\,\,\,\, t_0 < t < 2t_0
 ,
 $$
    one has
 $$
 F_1(t) > \frac{1}{A-1}+ hA^n,\,\,\, \text{for} \,\, t \,\, \text{from the interval}\,\,\, 2^n t_0 < t < 2^{n+1}t_0,
 $$
 and this means that 
 $$
 \lim_{t\to +\infty} 
  F_1(t) = +\infty
 .
 $$
 
    \vskip+0.3cm
    
    As for series
 $F_2(t)$,
 it satisfies inequalities
 $$
 -AF_2(t) - 1 < F_2(2t)< -AF_2(t) +1,\,\,\,\,\,\forall \, t.
 $$
  So if for certain $t, h >0$ one has
 $$
 F_2(t) >\frac{1}{A-1}+h
 ,$$
 then
 $$
 F_2(2^nt) > \frac{1}{A-1}+ hA^n,\,\,\,\, \,
 \,\,\,\text{for even}\, n
 $$
 and
  $$
 F_2(2^nt) < -\frac{1}{A-1} - hA^n,\,\,\,   
 \,\,\,\text{for odd}\, n.
 $$
 This means that 
  $$
 \liminf_{t\to +\infty} 
  F_2(t) = -\infty,\,\,\,\,\,\,\,\,
  \limsup_{t\to +\infty} 
  F_2(t) = +\infty
 .
 $$
   \vskip+0.3cm
   The results described above seem to have no relations to Diophantine Approximation and Kronecker sequence. However this relation exists.
   In order to explain this, we at first describe  an example form the next paper by  Poincaré \cite{P2} and then give an interpretation of the construction in terms of sums with Kronecker sequence.
 
    \vskip+0.3cm
 In \cite{P2}  Poincaré considered a similar series with coefficients of Diophantine nature.
 Let  $D$ be a squarefree positive integer and integers $u,v$ be solutions of Pell's equation
 $$
 u^2-Dv^2 = 1.
 $$
 Let 
 $$ \lambda = u-v\sqrt{D} \in (0,1), \,\,\,\,\,\,\, \lambda^n = u_n-v_n\sqrt{D},\,\,\, u_n, v_n \in \mathbb{Z}
 $$
 and consider the series
 \begin{equation}\label{ser}
 F_3(t) = 
 \sum_{n=1}^\infty A^n\sin \lambda^n t.
 \end{equation}
 Then according to the argument behind,   for $|A|>1 ,|\lambda A|<1$ series (\ref{ser}) converges non-uniformly.  Moreover, it   tends to infinity as $t\to \infty$ in the case when $A\in \left(1,\frac{1}{\lambda}\right)$ and
 { is close to $\frac{1}{\lambda}$}. In the case when $A\in  \left(-\frac{1}{\lambda},-1\right)$
 { is close to 
 $-\frac{1}{\lambda}$}
 the series (\ref{ser}) is unbounded and oscillates.  Poincaré obtained series (\ref{ser}) as a solution of a certain simple differential equation, in particular, it can be represented as an integral
 \begin{equation}\label{F3}
 F_3(2\pi t) = 2 \pi \int_0^t g(\omega_1s, \omega_2s) ds,
 \end{equation}
 {where}
 $
 \omega_1 = 1,\,\omega_2= \sqrt{D}$
 and
 \begin{equation}\label{gF3}
 g(x_1,x_2) = \sum_{n=1}^\infty (A\lambda)^n \cos(2\pi  (u_nx_1- v_nx_2)) 
\end{equation}
 is a continuous periodic function.
   \vskip+0.3cm

   Now we describe a discrete version of the last example. Instead of considering continuous  periodic function  $
 g(x_1,x_2) $ in two variables $x_1$ and  $x_2$  we  deal with $1$-periodic  continuous functions
 $$
 f(x) = \int_0^1 g(s, s\sqrt{D}+x) ds 
 $$
 in one variable $x\in \mathbb{R}$, where   $
 g(x_1,x_2) $
 is defined in (\ref{gF3}). Then by (\ref{F3}) for integer values of $t$ one has 
 $$
 \frac{
F_3(2\pi t)}{2\pi}  = 
\int_0^t g(s, s\sqrt{D}) ds=
\sum_{k=1}^t  \int_{k-1}^k g(s, s\sqrt{D}) ds
 = \sum_{k=0}^{t-1}   f(k\sqrt{D}).
$$
We can take $\alpha = \sqrt{D}$ with squarefree $D\in \mathbb{Z}_+$ and consider infinite Kronecker sequence 
\begin{equation}\label{kkr}
\{ k\alpha\}, \,\,\,\, k =1,2,3,.... .
\end{equation}
We see that
  Poincaré was the first to construct a continuous $1$-periodic function $f(x)$ with zero mean value
  $
  \int_0^1 f(t) dt = 0$ such that  for the sum of the values of $f$ over the sequence (\ref{kkr}) one has
  \begin{equation}\label{infi}
  \sum_{k=0}^{t-1}   f(k\alpha) \to +\infty,\,\,\,\,  t \to +\infty,
  \end{equation}
  and the proof of this fact relies on the argument described in the beginning of this section.

    \vskip+0.3cm

Poincar\'{e}'s
examples described  here are discussed in details in \cite{Kz}, Chapter 8. In the next section we will briefly discuss results related to the sequence (\ref{kkr}), its multidimensional analogs and corresponding sums.

    \vskip+0.3cm
   
    {\bf 2. On Birkhoff sums.}
    
       \vskip+0.3cm
    
    Let $\pmb{x} = (x_1,...,x_n) \in \mathbb{R}^n$ and $ f(\pmb{x})$ be a real-valued function which is 1-periodic in each variable $x_j$, that is $f:\mathbb{T}^n \to \mathbb{R}$ and
    $\mathbb{T}^n$ is a standard $n$-dimensional torus.  Let $\pmb{\alpha} = (\alpha_1,...,\alpha_n)\in \mathbb{R}^n$. The main object of the present paper is the sum 
\begin{equation}\label{summ0}
S^t_f (\pmb{\alpha},\pmb{x}) = \sum_{k=0}^{t-1} f(k\pmb{\alpha}+\pmb{x}).
\end{equation}
For the particular case $\pmb{x} = \pmb{0}$ we use the notation 
\begin{equation}\label{summ}
S^t_f (\pmb{\alpha}) = S^t_f (\pmb{\alpha},\pmb{0}) =  \sum_{k=0}^{t-1} f(k\pmb{\alpha}).
\end{equation}
We will define vector $\pmb{\alpha} = (\alpha_1,...,\alpha_n)$ to be {\it completely irrational} if
$1, \alpha_1,...,\alpha_n$ are linearly independent over $\mathbb{Q}$. 
We would like to give here several remarks and comments concerning these objects.
Some of the results mentioned  below in this section as well as some other results are discussed in \cite{K} and in  the last parts of  survey \cite{mou}.

   \vskip+0.3cm
 {\bf A.}
According to the fundamental H. Weyl's theorem for Riemann integrable $ f:\mathbb{T}^n \to \mathbb{R}$ 
in the case when  $\pmb{\alpha}\in \mathbb{R}^n$  is completely irrational, for any $\pmb{x}$  we have asymptotic equality
$$
\frac{1}{t}\, S^t_f (\pmb{\alpha},\pmb{x}) =  \frac{1}{t}\,\sum_{k=0}^{t-1} f(k\pmb{\alpha}+\pmb{x}) \to  \int_{\mathbb{T}^n} f(\pmb{z})d \pmb{z}\,\,\,\,\,\
\text{when}\,\,\,\,\, t \to \infty.
$$

   \vskip+0.3cm
 {\bf B.} 
 A wonderful result was obtained by Marstrand \cite{mr}. For any $\varepsilon >0$ there exists  an open set $\Omega \subset [0,1]$  of Lebesgue measure 
 $\mu \Omega <\varepsilon$ 
 such that for 1-periodic function $f_\Omega (t) $ which coincide with characteristic function of the set $\Omega$ on the interval $[0,1]$ and
 {\it any} real $\alpha$ one has
 $$
 \limsup_{t\to \infty} \,
 \frac{1}{t}\, \sum_{k=1}^t   f_\Omega(k{\alpha})  =1 \ne \int_0^1f_\Omega (t) dt = \mu \Omega.
 $$
 We see that in this construction $f_\Omega (t) \in L^1 [0,1] $ but is not a continuous function and even
 Riemann integrable function.
 
   \vskip+0.3cm
 {\bf C.} 
 For $n=1$ and a smooth function $f$ for irrational $\alpha$  the limit behaviour (\ref{infi}) never holds.
 If $f$ is a  1-periodic function of bounded variation ${\rm V}[f]$ and $\int_0^1 f(x)dx=0$, Koksma's inequality (see \cite{KN}, Ch.2, \S 5)  ensures the bound
 $$
 |S^t_f (\alpha) |\le {\rm V}[f] D_t,
 $$
 where 
 $$
 D_t =
 \sup_{\gamma\in [0,1]} |\#\{ k \in \mathbb{Z}_+: \, 1\le k \le t,\,\,\{k\alpha\}< \gamma\} - t\gamma|
 $$ is the discrepancy of the sequence $ \{k\alpha\} , 1\le k\le t$. When $ t = q_n$ is a denominator of a convergent to $\alpha$, one has
 $ D_{q_{n}} =O(1).$ So for  function $f$ of bounded variation and irrational $\alpha$ we have
 $$
 \liminf_{t\to \infty}  \max_{x\in [0,1]} |S^t_f (\alpha, x) |<+\infty.
 $$
 If we consider an absolutely continuous function  (this condition is stronger the the condition of bounded variation, see \cite{KF}, Ch. 9), we have a result by Sidorov \cite{sy}:
 for any  1-periodic absolutely continuous function $f$  with  zero mean value and any irrational $\alpha$ one has
  $$
 \liminf_{t\to \infty} \max_{x\in [0,1]} |S^t _f(\alpha, x) | =0.
 $$
 There is no reasonable analog of Sidorov's result for $n\ge 2$, see the discussion in the beginning of Subsection 8.2 below.
 %Earlier results in this direction are due to Kozlov \cite{ko} and  Krygin \cite{kry}.

    \vskip+0.3cm
    
     {\bf D.} In the case $n=1$ there are many generalisations of Poincare's example  which deal with the rate of growth of  the sums (\ref{summ0}),
    various conditions on $f$, Hausdorff dimension of the set of initial phases $x$ for which one has growth of sums (\ref{summ0})  (see \cite{bes,12,11,13,14}).
    For example, Besicovich \cite{bes} showed that  for any irrational $\alpha$ there exists continuous 1-periodic $f$ with zero mean value, such that 
    $$
     \liminf_{t\to \infty} \frac{S^t_f (\pmb{\alpha})}{\sqrt{t}} >0. 
    $$
    In some sense, the strongest result for continuous  $f$ is due to Kochergin \cite{K}. We will formulate it in the next section.

          \vskip+0.3cm
    
     {\bf E.}  Recently Colzani \cite{Co} found very natural  sufficient conditions on periodic $f:\mathbb{R}^n \to \mathbb{R} $ with zero mean value to ensure  that 
     $$
     \sup_{t\in \mathbb{Z}_+}\, 
     \max_{x\in [0,1]} |S^t _f(\alpha, x) |  <\infty
     $$
     for almost all $\pmb{\alpha}\in \mathbb{R}^n$. We wold like to formulate  some of his results in a separate Section 4 below.
     
    \vskip+0.3cm
    
     {\bf F.} Let $w\in \mathbb{Z}_+$.
     We define function $g(\pmb{x}) = g(x_1,...,x_n):\mathbb{R}^n\to \mathbb{R}$ to be {\it $w$-smooth}
     if all its partial derivatives
     $$
     \frac{\partial g^w}{\partial^{w_1}x_1\,\,  ...  \,\, \partial^{w_n}x_n},\,\,\,\,\,\,\, 0\le w_1,...,w_n\le w,\,\,w_1+...+w_n = w
     $$
     are continuous functions in $\mathbb{R}^n$.
     
     In \cite{Mnon} the author for any  badly approximable $\pmb{\alpha}$ (see definition below in Section 4)   constructed a $(n-1)$-smooth periodic  function $f(\pmb{x}):\mathbb{R}^n\to \mathbb{R}$ with zero mean value such that  
     \begin{equation}\label{lininf}
     \inf_{t>0} \, S^t_f (\pmb{\alpha}) >0.
     \end{equation}
     Later Konyagin \cite{kon} constructed examples of completely irrational  $\pmb{\alpha}\in \mathbb{R}^n$ and 
     $w$-smooth function and $\pmb{\alpha}\in \mathbb{R}^n$  with 
     $w =w_n =\left[\frac{2^n(n-1)^{n-1}}{n^n}\right]-1\sim  \frac{2^n}{en}, \, n\to \infty $ such that  (\ref{lininf}) holds
 (both papers \cite{Mnon, kon} deal with integrals, so the number of variables in these paper is equal  to $n+1$ in our notation). We should note here that 
    $w_n> n-1$ for $n\ge 8$. In the present paper we use bounds with Diophantine exponents to deduce these results and even to give a small improvement.
    In Section 8 we give an example of $\pmb{\alpha}\in \mathbb{R}^n$  and $w$-smooth function in $n$ variables  with  condition  (\ref{lininf}) and
     $ w\sim \frac{2^{n+1}}{en}, \, n\to \infty $.

     It is known (see \cite{m1,m2})  that if 
   function $f(\pmb{x})$ in $n$ variables with zero mean value is  $w$-smooth with $ w\asymp n^{ cn}$ with large $c$, then for any completely irrational $\pmb{\alpha}$ one has
     $$
      \liminf_{t\to \infty} \, |S^t_f (\pmb{\alpha})| <\infty.
     $$
     We do not discuss here this complicated result and refer to  the last part of survey \cite{mou}.
     
    \vskip+0.3cm

     {\bf G.} Sums (\ref{summ0}) are particular cases of Birkhoff sums for the uniquely ergodic rotation of torus $ \pmb{x}\mapsto\pmb{x}+\pmb{\alpha}$.
     There are some general theorems
     \cite{Xale,Halas,Peres,Krygin,Rysikov} 
      concerning almost all trajectories of dynamical systems which can give  a lot of information about  behaviour of sums (\ref{summ0})
     for generic $\pmb{x}$.
    \vskip+0.3cm

    {\bf 3. A result by Kochergin.}

   \vskip+0.3cm

   Here we formulate the main result from a recent paper \cite{K}. Let $\alpha \in \mathbb{R}\setminus\mathbb{Q}$ and
   $ \Sigma = \{\sigma_t\}_{ t\in \mathbb{Z}_+}$ be a sequence of real numbers decreasing to  zero. Then there exists a continuous
   1-periodic function  $f:\mathbb{T}\to\mathbb{R}$ with  
   $$
   \int_0^1 f(x) dx = 0,\,\,\,\,\,\,\, \max_{x\in [0,1]} |f(x)| = 1,
   $$
    a set $\frak{D}\subset [0,1]$ of Hausdorff dimension
   ${\rm dim}_{\rm H} \,\frak{D} = 1$ and $t_0$, such that 
$$
S_f^t (\alpha, x) \ge t\sigma_t\,\,\,\,\,\,\, \,\,\,\forall \, t\ge t_0,\,\, \forall\, x\in \frak{D}.
$$
Of course, here $f, \frak{D}$ and $t_0$ depend on $\alpha $ and $\Sigma$.

   \vskip+0.3cm
   
A similar result  without the condition on thickness of the set $\frak{D}$ (just an existence result) was published a little bit earlier in \cite{A}.
In the present paper we are mostly interested just in existence-type results.

   \vskip+0.3cm
  In papers \cite{A,K} function $f$ is constructed   as a sum of a series which converges uniformly and which elements are piecewise linear functions.

   \vskip+0.3cm

      \vskip+0.3cm
      
    {\bf 4. Results by Colzani.}
    
          \vskip+0.3cm
          We would like to formulate here certain  results
          by Colzani 
          from a recent  paper \cite{Co}.

       For vectors 
      $
\pmb{x} = (x_1,...,x_n),
\,\,
\pmb{m} = (m_1,...,m_n) \in \mathbb{R}^n$
we use the notation
$$
\pmb{x}.\pmb{m} = x_1m_1+...+x_nm_n.
$$
For  periodic function $ f:\mathbb{R}^n \to \mathbb{R}$ with zero mean value we consider its Fourier series
\begin{equation}\label{fgf}
f(\pmb{x}) \sim \sum_{\pmb{m}\in \mathbb{Z}^n\setminus\{\pmb{0}\}}
\,\,
f_{\pmb{m}} \exp (2\pi i  \pmb{m}.\pmb{x})
\end{equation}
    
          \vskip+0.3cm
The  main result from \cite{Co} is as follows. Consider positive function $\Phi(t):[1,\infty)\to \mathbb{R}$ such that  both functions $ t\mapsto \Phi(t)$ and $t\mapsto 1/\Phi(t)$ are increasing and
\begin{equation}\label{io}
\int_2^{+\infty} \frac{dt}{t\Phi(t) }<+\infty
\end{equation}
and assume that Fourier series (\ref{fgf}) converges absolutely,
that is
\begin{equation}\label{ggggg}
 \sum_{\pmb{m}\in \mathbb{Z}^n\setminus\{\pmb{0}\}}
 |
f_{\pmb{m}}
|
<+\infty.  
\end{equation}
Then for almost all $\pmb{\alpha}\in \mathbb{R}^n$ (in the sense of Lebesgue measure) 
there exists positive $c(f,\pmb{\alpha})$ such that 
for every $\pmb{x}\in \mathbb{R}^n$ one has
$$
\left|S_f^t (\alpha, x) \right| \le c(f,\pmb{\alpha}) \Phi (t)\,\,\,\,\, \forall \, t \ge 1.
$$

          \vskip+0.3cm
          For example, to satisfy (\ref{io}) one can take
          $\Phi (t) = \log(t+1) \log^{1+\varepsilon}(\log t+2)) $ for arbitrary $\varepsilon >0$.
          
           \vskip+0.3cm
           In the same paper under a little bit more restrictive condition
           $$
            \sum_{\pmb{m}\in \mathbb{Z}^n\setminus\{\pmb{0}\}}
 |
f_{\pmb{m}}
| \log \left(1+\frac{1}{
|f_{\pmb{m}}|}
\right)
<+\infty
$$
Colzani proves  that  the sum
$
S_f^t (\alpha, x) 
$
is bounded as $ t\to \infty$ for every $\pmb{x}$ and for almost all $\pmb{\alpha}$.

      \vskip+0.3cm
      At the end of this section we would like to formulate one more result from the paper \cite{Co} which deals just with the one-dimensional case . It shows that  condition (\ref{ggggg}) on absolute convergence is close to an optimal one. Suppose that 1-periodic function  $f(x):\mathbb{R}\to\mathbb{R}$ with zero mean value is square integrable and the sequences of its Fourier coefficients $|f_m|, m =1,2,3,...$ and
      $|f_{-m}|, m =1,2,3,...$  are decreasing. Then under the condition that there exists a set of positive Lebesgue measure $A \subset [0,1]$ such that 
      $$
      \int_0^1 |S_f^t(\alpha, x)|^2 dx <+\infty\,\,\,\,\, \forall \alpha\in A
      ,$$
       function $f$ has absolutely convergent Fourier expansion, that is (\ref{ggggg}) holds.

.

   \vskip+0.3cm
   
    {\bf 5. Diophantine approximation.}

   \vskip+0.3cm

    For completely irrational  $\pmb{\alpha} = (\alpha_1,...,\alpha_n)\in \mathbb{R}^n$ we consider irrationality measure function
   $$
   \psi_{\pmb{\alpha}} (t ) = \min_{\pmb{m}\in \mathbb{Z}^n:\, 0< \max_j|m_j|\le t}\,\,\,\, || \pmb{\alpha}.\pmb{m} ||,\,\,\,\,\,
   \text{where}\,\,\,\,\, ||\cdot || = \min_{a\in \mathbb{Z}}|\cdot - a|.
   $$
   Irrationality measure function is a piecewise constant function which defines best approximation vectors
   $$
   \pmb{m}_\nu = (m_{1,\nu},...,m_{n,\nu}) \in \mathbb{Z}^n,\,\,\,\, M_\nu =     \max_{1\le j \le n}|{m}_{j,\nu} | 
   $$
   such that 
   $$
      \psi_{\pmb{\alpha}} (t ) = || \pmb{\alpha}.\pmb{m}_\nu || = :\zeta_\nu\,\,\,\, \text{for}
      \,\,\,\,
M_\nu\le t < M_{\nu+1}.
      $$
      All the general facts about the behaviour of  best approximation vectors are discussed in   \cite{che}
      and \cite{mou}.
      In particular, it follows from the Minkowski convex body theorem that
   \begin{equation}\label{mink}
   \zeta_\nu \le  M_{\nu+1}^{-n}.
   \end{equation}
   It is important for us that $M_\nu$ grow exponentially and so
   \begin{equation}\label{serri}
   \sum_{\nu=1}^\infty M_\nu^{-1} < \infty.
   \end{equation}
   An important tool to estimate the quality of the best approximation is related to the ordinary and the uniform Diophantine exponents.
   The {\it ordinary Diophantine  exponent}  $\omega = \omega(\pmb{\alpha})$ is defined by
   $$
   \omega = \omega(\pmb{\alpha})= \sup\{ \gamma \in \mathbb{R}: \, \liminf_{t\to \infty}  t^\gamma \cdot \psi_{\pmb{\alpha}} (t )<\infty \},
   $$
   while the  {\it uniform Diophantine exponent}
    $\hat{\omega} = \hat{\omega}(\pmb{\alpha})$ is defined by
   $$
   \hat{
   \omega} = \hat{\omega}(\pmb{\alpha})= \sup\{ \gamma \in \mathbb{R}: \, \limsup_{t\to \infty}  t^\gamma \cdot \psi_{\pmb{\alpha}} (t )<\infty\}.
   $$
   From the definitions and (\ref{mink}) we see that 
   $$
   n\le \hat{\omega}\le \omega.
   $$
   \vskip+0.3cm
     In terms of Diophantine exponents $\omega,\hat{\omega}$ we can formulate upper and lower bounds for the values  $\zeta_\nu$ of the best approximations.
   Namely, 
     \begin{equation}\label{be1}
     \forall\, \delta >0 \,\,\,\,\, \forall \, \nu \,\text {lagre enough}\,\,\,\,\,\,\,\text{one has}\,\,\,   M_\nu^{-\omega-\delta}\le  \zeta_\nu \le M_{\nu+1}^{-\hat{\omega}+\delta}.
     \end{equation}

      \vskip+0.3cm

   A result by Marnat and Moshchevitin \cite{mamo} gives an optimal bound for the ratio of the exponents under consideration. Suppose that $\pmb{\alpha}$ is completely irrational. Then
     \begin{equation}\label{mamo}
   \frac{\omega}{\hat{\omega}} \ge G (\hat{\omega}),
   \end{equation}
   where  $ G (\hat{\omega}) \ge 1$ is the unique positive root of the polynomial
   \begin{equation}\label{mamo1}
   R_{n,\hat{\omega}}^*(x) = x^{n-1}+x^{n-2}+...+ x + 1 -\hat{\omega}
   .
   \end{equation}
   Optimality of the inequality (\ref{mamo}) means following. For any real numbers $ t, \tau$ satisfying inequalities
   $$
   \tau\ge n,\,\,\,\, t \ge G(\tau)
   $$
   there exists completelly irrational $\pmb{\alpha}\in \mathbb{R}^n$ with
   $$
   \hat{\omega}(\pmb{\alpha}) = \tau,\,\,\,\, 
      \frac{\omega (\pmb{\alpha})}{\hat{\omega}(\pmb{\alpha})}  = t.
      $$
      In Section 8 below we should consider  quantities $ 2\hat{\omega}(\pmb{\alpha})- {\omega}(\pmb{\alpha})$ and
       $$
       \frak{g}_n := \max_{\pmb{\alpha}\in \mathbb{R}^n -\,\,\text{completely irrational}}
       \,\,\,\,\,\,\,
       (2\hat{\omega}(\pmb{\alpha})- {\omega}(\pmb{\alpha})).
       $$
     For this purpose we need the following
     
        \vskip+0.3cm

      {\bf Proposition 1.} \,\,\,\,\,$\frak{g}_n \ge  \frak{g}_n^*:=\frac{2}{n-1} \, \left(2^{n} \left( 1-\frac{1}{n+1}\right)^n - 1\right) \sim \frac{2^{n+1}}{en},\,\,\,n\to \infty$.

     \vskip+0.3cm
     
Proof. From the optimality of inequality (\ref{mamo}) and the explicit formula (\ref{mamo1}) for polynomial $    R_{n,\hat{\omega}}^*(x) $  it follows that 
$$
\frak{g}_n  =\max_{\tau \ge n}  \, \tau \cdot (2-G(\tau)) =
\max_{t\in (1,2)}\, (1+t+t^2 + ...+t^{n-1}) (2-t) =
\max_{t\in (1,2)}\, g(t) \ge g\left(2- \frac{2}{n+1}\right) = \frak{g}_n^*,
$$
where $g(t) = \frac{(t^n-1)(2-t)}{t-1}$.$\Box$

        \vskip+0.3cm
      We  should formulate the definition of uniform Diophantine exponent for simultaneous approximation which will be important in Theorem 5 in Section 8 below. Consider irrationality measure function for simultaneous approximation
     $$
     \psi^*_{\pmb{\alpha}}(t) = 
     \min_{q\in \mathbb{Z}_+:\, q\le t }\,\, \max_{1\le j \le n} ||q\alpha_j||.
     $$
     {\it Uniform Diophantine exponent for simultaneous approximation} is defined as
     $$
     \hat{\lambda} =\hat{\lambda}(\pmb{\alpha}) =
     \sup\{ \gamma \in \mathbb{R}: \, \limsup_{t\to \infty}  t^\gamma \cdot \psi_{\pmb{\alpha}}^* (t )<\infty\}.
   $$
   Many properties of $\hat{\lambda}(\pmb{\alpha}) $ are discussed in  \cite{mamo,mou}.
   We do not want go into further details here. We only should note that  for $\pmb{\alpha}\not \in \mathbb{Q}^n$  the inequalities
   $$
   \frac{1}{n} \le \hat{\lambda}(\pmb{\alpha})\le 1
   $$
   are valid, and that optimal transference inequality
   \begin{equation}\label{transfer}
   \hat{\lambda}(\pmb{\alpha})\ge \frac{1- \frac{1}{\hat{\omega}(\pmb{\alpha})}}{n-1}
   \end{equation}
  between $\hat{\omega}(\pmb{\alpha})$ and $\hat{\lambda}(\pmb{\alpha})$ for $ n \ge2$
   was obtained by German in \cite{g1,g2}, meanwhile for $n=2$ we have a famous equality
   $$
  \hat{\lambda}(\pmb{\alpha}) = 1- \frac{1}{\hat{\omega} (\pmb{\alpha})}
   $$
   due to  Jarn\'{\i}k \cite{j}. Optimality of (\ref{transfer}) was established by Marnat in \cite{Mar} and independently by Schmidt and Summerer \cite{SSS}. For further information see recent survey \cite{GGG}.

     \vskip+0.3cm
     At the end of this section we would like to recall the reader that $\pmb{\alpha}$ is called {\it badly approximable} if
     $$
          \inf_{t\ge 1} \,\,t^{n} \cdot \psi_{\pmb{\alpha}	} (t) >0,
$$
or equivalently
     \begin{equation}\label{badd}
     c(\pmb{\alpha}) M_\nu^{-n} \le 
    \zeta_\nu \le M_\nu^{-n} \,\,\,\, \forall\, \nu 
     \end{equation}
     with some positive $     c(\pmb{\alpha})$.
     The set of all badly approximable vectors in $\mathbb{R}^n$ has zero Lebesgue measure in $\mathbb{R}^n$ but full Hausdorff dimension
     (see \cite{sch}, Chapter 3). For badly approximable $\pmb{\alpha}\in \mathbb{R}^n$ one has 
     $\omega(\pmb{\alpha}) = \hat{\omega}(\pmb{\alpha}) = n$ and $\hat{\lambda}(\pmb{\alpha}) =1/n$.

    \newpage
   
    {\bf 6. Multidimensional version of Kochergin's result.}

   \vskip+0.3cm
   Here we give a multidimensional version of Kochergin's theorem. We are interested here only in the behaviour of a particular sum with $ \pmb{x} = \pmb{0}$ for the rotation of torus $\mathbb{T}^n$.
   Our proof follows the main argument from \cite{A,K} and  is related to the best  Diophantine approximations.
 
   \vskip+0.3cm

{\bf Theorem 1.} {\it
For any  completely irrational vector  $\pmb{\alpha}$  and for any sequence  $ \sigma_r \downarrow 0$ there exists a continuous function with zero mean value
$f(\pmb{x}) :\mathbb{T}^n\to \mathbb{R}$ and $t_0 \in \mathbb{Z}_+$  such that
\begin{equation}\label{po}
\max_{\pmb{x}\in \mathbb{T}^n} |f(\pmb{x})| \le \pi.
\end{equation}
and
\begin{equation}\label{iu}
S_f^t (\pmb{\alpha}) =
\sum_{k=0}^{t-1} f(k\pmb{\alpha }) 
 \ge
 t\sigma_t\,\,\,\,\,\,\forall\, t\ge t_0.
\end{equation}
}

\vskip+0.3cm

Proof.
 We use best approximation vectors to $\pmb{\alpha}$  to define  functions
\begin{equation}\label{funkktion}
F_\nu (\pmb{x}) =   
\frac{|\sin\, \pi (\pmb{m}_\nu.\pmb{x}) | }{ \zeta_\nu},
\end{equation}
and
$$
f_\nu(\pmb{x}) = F_\nu(\pmb{x}+\pmb{\alpha}) - F_\nu(\pmb{x})
=
 \frac{|\sin\, \pi (\pmb{m}_\nu.(\pmb{x} + \pmb{\alpha})) | - |\sin \pi (\pmb{m}_\nu.\pmb{x})|}{ \zeta_\nu}=
\frac{|\sin\, \pi (\pmb{m}_\nu.\pmb{x} + \zeta_\nu) | - |\sin \pi (\pmb{m}_\nu.\pmb{x})|}{ \zeta_\nu}.
$$
As
$$
\frac{\left|\,|\sin\, \pi (\pmb{m}_\nu.\pmb{x} + z) | - |\sin \pi (\pmb{m}_\nu.\pmb{x})| \,\right|}{ |z|}
\le \pi
,$$
we see that 
\begin{equation}\label{pi}
 \max_{\pmb{x}\in \mathbb{T}^n} |f_\nu(\pmb{x})| 
\le \pi
\end{equation}
and
\begin{equation}\label{pi1}
 \sum_{k=0}^{t-1} f_\nu(k\pmb{\alpha}) = F_\nu (t\pmb{\alpha}).
\end{equation}

We define sequence $r_k, k=1,2,3,...$ by the following procedure. For a fixed $k$ we take large $\nu_k$ to satisfy
\begin{equation}\label{2_jj}
\sigma_{r_k}\le \sum_{j=k+1}^\infty 2^{-j},\,\,\,\,\, r_k = \left[\frac{1}{2\zeta_{\nu_k}}\right].
\end{equation}
Also we assume than in our construction $r_k < r_{k+1}$. We can ensure this inequality while choosing $r_{k+1}$ after $r_k$.
Now we define functions
\begin{equation}\label{2_j}
F(\pmb{x})  = \sum_{j=1}^\infty 2^{-j} F_{\nu_j}  (\pmb{x}) ,\,\,\,\,\,
f(\pmb{x})  = \sum_{j=1}^\infty 2^{-j} f_{\nu_j}  (\pmb{x}) .
\end{equation}
(Here we should note that instead of coefficients $2^{-j}$  in (\ref{2_j})  and in upper bound  (\ref{2_jj})  for $\sigma_{r_k}$ one can take an arbitrary sequence   of positive integers $ a_j>0$ satisfying $\sum_{j=1}^\infty a_j = 1$ as it was done in \cite{K}.)
By (\ref{pi}), the second series here  converges uniformly and
$
f(\pmb{x})$ is a continuous function with zero mean value, so 
  we have (\ref{po}).
By (\ref{pi1}) we get
\begin{equation}\label{sss}
\sum_{k=0}^{t-1} f(k\pmb{\alpha }) = F(t\pmb{\alpha }). 
\end{equation}

 Now we explain why inequality (\ref{iu}) is valid. For integer $t$ from the interval   $r_k \le t < r_{k+1}$ we have $ t< \frac{1}{2\zeta_{\nu_{k+1}}}$ and so  as $ F_\nu (\pmb{x}) \ge 0$ and $
|\sin(\pi z)| \ge 2 ||z||
$
we get
$$
F(t\pmb{\alpha}) \ge \sum_{j=k+1}^\infty 2^{-j}   F_{\nu_j}  (t\pmb{\alpha})=
\sum_{j=k+1}^\infty 2^{-j}  \frac{|\sin (\pi t\zeta_{\nu_j})|}{\zeta_{\nu_j}}\ge
2t\sum_{j=k+1}^\infty 2^{-j} \ge 2t\sigma_{r_k} \ge2t\sigma_t.
$$
The last inequality together with (\ref{sss})  gives us (\ref{iu}). Theorem is proven.$\Box$

     \vskip+0.3cm
   
    {\bf 7. Nonexistence of "universal"  continuous function $f$.}

   \vskip+0.3cm
   
   For  a continuous function its Fourier series may not converge absolutely. So we cannot apply to continuous functions the result by Colzani from section 4 .
   Nevertheless,
   here we show that it is not possible to find a continuous function $f$  which will be the same for all irrational $\pmb{\alpha}$. We formulate here a very easy result which deals with  arbitrary value of $n$ (Theorem 2)  and a slightly better result for $n=1$ (Theorem 3).
   \vskip+0.3cm

{\bf 
Theorem 2.} {\it  Let $ n \ge 1$.
Let $ \psi(t) \to +\infty$ as $t\to \infty$.
For any continuous $1$-periodic (in each variable) function $f(\pmb{x}):\mathbb{R}^n\to \mathbb{R}$  with zero mean value
\begin{equation}\label{II1}
\int_{\mathbb{T}^n} f(\pmb{x}) d\pmb{x} = 0
\end{equation}
there exists  a completely irrational 
$\pmb{\alpha} \in \mathbb{R}^n$ such that  
\begin{equation}\label{II2}
\liminf_{t\to \infty} \frac{S^t_f(\pmb{\alpha})}{\psi(t)\cdot \log t} = 0.
\end{equation}
}
 
%\vskip+0.3cm
%{\bf Remark.}  We believe this statement to be true for constant function $ \psi (t)$ but we cannot prove it now even with $<+\infty$ instead of $0$ in r.h.s.
%of (\ref{II2}).
%\vskip+0.3cm

For $n= 1$ we have a better result.

\vskip+0.3cm

{\bf 
Theorem 3.} {\it 
Let $ \psi(t) \to +\infty$ as $t\to \infty$.
For any continuous $1$-periodic function $f:\mathbb{R}\to \mathbb{R}$  with zero mean value
there exists 
$\alpha \in \mathbb{R}\setminus \mathbb{Q}$ such that  
\begin{equation}\label{II2}
\liminf_{t\to \infty} \frac{S_f^t(\alpha)}{\psi(t)} = 0.
\end{equation}
}

   \vskip+0.3cm

In the proofs below it is more convenient to consider the sum
$$
S^t(\pmb{\alpha}) =
S^t_{f_*}(\pmb{\alpha}) =  \sum_{k=1}^t f (k\pmb{\alpha})
$$
rather than the sum (\ref{summ}).  This will cause no problem.
   \vskip+0.3cm

Proof of Theorem 2.
Define
$$
M = \max_{\pmb{x}\in \mathbb{T}^n} 
|f(\pmb{x})|.
$$
We consider a box
$$
J = [\beta_1,\gamma_1]\times ...\times  [\beta_1,\gamma_1]\subset [0,1)^n,\,\,\,\,\, \beta_j <\gamma_j
$$
 and an integral
$$
I=\int_{J}S^t (\pmb{x}) d\pmb{x} =
\sum_{k=1}^t
\int_{J}f (k\pmb{x})  d\pmb{x}=
\sum_{k=1}^t\frac{1}{k^n}\, 
\int_{kJ}f (\pmb{x})  d\pmb{x},
\,\,\,\,\,
kJ = [k\beta_1,k\gamma_1]\times ...\times  [k\beta_1,k\gamma_1].
$$
Define $ \lambda_j = [k(\gamma_j -
\beta_j]$. Then by (\ref{II1}) we have
$$
\int_{kJ}f (\pmb{x})  d\pmb{x} = \int_{J_*}f (\pmb{x})  d\pmb{x},\,\,\,\, J_* = kJ \setminus [k\alpha_1, k\alpha_1 + \lambda_1]\times ...\times  [k\alpha_n, k\alpha_n + \lambda_n].
$$
It is clear that 
$$
J_* = \cup_{i=1}^n J_* (i),\,\,\,\,\text{where}\,\,\,\,\,
J_* (i) =\{ (x_1,...,x_n)\in \mathbb{R}^n:  k\alpha_i + \lambda_i \le x_i \le k\beta_i,\,\,\,  k\alpha_l \le x_l \le k\beta_l,\,\, l \neq i\}.
$$
So
$$
\left|
\int_{J_*}f (\pmb{x})  d\pmb{x}
\right|
\le M \int_{J_*}  d\pmb{x} = M \mu J_*\le
M\sum_{i=1}^n \mu J_* (i) \le Mn
k^{n-1},
$$
and
\begin{equation}\label{log}
\left|\int_{J}S^t (\pmb{x}) d\pmb{x} \right| \le Mn(1+ \log t).
\end{equation}

Now we construct a sequence of integers $t_\nu$
and a sequence of nested boxes
$$
J_\nu \supset J_{\nu+1}
$$
such that 
\begin{equation}\label{tn1}
|S^{t_{\nu-1}}(\pmb{x})| \le \frac{2Mn}{\mu J_{\nu-1}} \,(1+ \log t_{\nu-1})\,\,\,\,\,\,\,\, \forall\, \pmb{x} \in J_{\nu}
\end{equation}
and
\begin{equation}\label{tn2}
 {\mu J_{\nu}} \cdot\psi (t_\nu)\to \infty,\,\,\,\,\, t\to \infty.
 \end{equation}
 It is clear that we can choose the boxes $J_\nu$ to avoid all the rational affine subspaces in $\mathbb{R}^n$.
This will prove Theorem 2  because  we get completely irrational  $\pmb{\alpha} \in \cap_\nu J_\nu $, and  condition (\ref{II2}) follows from (\ref{tn1}) and (\ref{tn2}).

We proceed the construction by induction. We begin  with $\nu=1$ and $ t_0 = 1$, $ J_0= J_1= [0,1]^s$. 
Suppose that  
$J_1,...,J_{\nu-1}, J_\nu$ and $t_0,t_1,...,t_{\nu-1}$ are constructed. We explain how to construct $ J_{\nu+1}$ and $t_{\nu}$.
By (\ref{log}),
$$
\left|\int_{J_{\nu}}S^t (\pmb{x}) d\pmb{x} \right| \le Mn\, (1+\log t)
$$
for all $t>1$.
If $S^{t}(\pmb{x})$ does not change sign on $J_{\nu}$, then
$$
\mu J_{\nu} \cdot \min_{\pmb{x}\in J_{\nu}}
|S^{t}(\pmb{x})| \le \left|\int_{J_{\nu}}S^t (\pmb{x}) d\pmb{x} \right| \le Mn \,(1+ \log t).
$$ 
If $S^{t}(\pmb{x})$ changes its  sign on $J_{\nu}$, then by continuity, $S^{t}(\pmb{x}_*)=0$ for some $ \pmb{x}_* \in J_{\nu}$.
In both cases for any $t>1$ we found  $ \pmb{x}_*(t) \in J_{\nu}$ such that 
$$
|S^{t}(\pmb{x}_*(t))| \le \frac{Mn}{\mu J_{\nu}}\,  (1+ \log t) .
$$
By continuity, for any $t$ exists a subbox $J_{\nu+1}= J_{\nu+1}(t)\subset J_\nu$, such that 
\begin{equation}\label{dzs}
\max_{\pmb{x}\in J_{\nu+1}}
|S^{t}(\pmb{x})| \le \frac{2Mn}{\mu J_{\nu}} \, (1+ \log t),
\end{equation}
and we see that 
(\ref{tn1}) holds for all $\pmb{x}\in J_{\nu+1}=J_{\nu+1}(t)$.
Box $J_\nu$ was defined at the previous step. So at the current step we can  choose $t_{\nu}=t$ large enough to satisfy (\ref{tn2}). 
Our inductive process is completed.

The inductive process described  proves Theorem 2.   
We can enforce $\pmb{\alpha}$ to be completely irrational by the following argument. Let us enumerate the countable set of all rational affine subspaces of $\mathbb{R}^n$ by
$\mathcal{A}_1,\mathcal{A}_2, \mathcal{A}_3,...$. Then at each step of the procedure we can choose $J_\nu$ to satisfy $ J_\nu\cap\mathcal{A}_\nu =\varnothing$. Then $ \pmb{\alpha}\not\in \cup_\nu \mathcal{A}_\nu$
It is clear that for  $\pmb{\alpha} \in \cap_\nu J_\nu $ condition (\ref{II2}) follows from (\ref{tn1}) and (\ref{tn2}).$\Box$

   \vskip+0.3cm

      \vskip+0.3cm

Proof of Theorem 3.

We will use the notation from the proof of Theorem 2.
Let
\begin{equation}\label{qII3}
a,q \in \mathbb{Z}_+, \,\,\,\, (a,q) = (a+1,q) = 1.
\end{equation}
We consider the integral
$$
I= \int_{\frac{a}{q}}^{\frac{a+1}{q}} S^t (x) dx =
\sum_{k=1}^t
\int_{\frac{a}{q}}^{\frac{a+1}{q}} f (kx) dx .
$$
By (\ref{II1}) we see that
$$
\int_{\frac{a}{q}}^{\frac{a+1}{q}} f (kx) dx  =
\frac{1}{k}\, 
\int_{\frac{a}{q}k}^{\frac{a+1}{q}k} f (x) dx =
\frac{1}{k}\, 
\int_{\left\{\frac{a}{q}k\right\}}^{\left\{\frac{a+1}{q}k\right\}} f (x) dx .
$$
So if we consider $1$-periodic function
$$
h(z) = \int_0^z f(x)dx,
$$
we have
$$
I = \sum_{k=1}^t \frac{w_k}{k},\,\,\,\,\,\,
w_k = 
\int_{\left\{\frac{a}{q}k\right\}}^{\left\{\frac{a+1}{q}k\right\}} f (x) dx = 
h\left( \frac{a+1}{q}k\right) - h\left( \frac{a}{q}k\right).
 $$
 It is clear that 
 \begin{equation}\label{qII4}
 |w_k| \le M.
 \end{equation}
 By partial summation we get
 \begin{equation}\label{qII6}
 I = \sum_{k=1}^t \left(\frac{1}{k}-\frac{1}{k+1}\right) W_k + \frac{W_t}{t+1},
 \,\,\,\,\,
 W_k = \sum_{j=1}^k w_j  .
 \end{equation}
 By periodicity of $h$  and the condition (\ref{qII4}) we see that 
 for the sum over all residues modulo $q$ we have
  $$
  \sum_{j\pmod{q}}
w_j
=  \sum_{j\pmod{q}} h\left( \frac{a+1}{q}j\right) -
 \sum_{j\pmod{q}} h\left( \frac{a}{q}j\right)
 = 0,\,\,\,\,\text{
and so}
\,\,\,\,\,
|W_k| \le M(q-1).
$$
From (\ref{qII6}) for any $t$ we get the bound
\begin{equation}\label{qII7}
|I| 
\le
 \sum_{k=1}^t \left(\frac{1}{k}-\frac{1}{k+1}\right)| W_k| + \frac{|W_t|}{t+1}
\le 
M(q-1) \cdot \left( \sum_{k=1}^t 
\left(\frac{1}{k}-\frac{1}{k+1}\right) + \frac{1}{t+1}
\right)
=  M(q-1).
\end{equation}
The rest of the proof is similar to those from Theorem 2.
By an inductive procedure we construct sequences
\begin{equation}\label{qII8}
q_1\le t_1<...<q_\nu<t_\nu<q_{\nu+1}<t_{\nu+1}<...
\end{equation}
and a sequence of nested segments
$$
J_\nu = \left[\frac{a_\nu}{q_\nu}, \frac{a_{\nu}+1}{q_{\nu}}\right]\supset J_{\nu+1}
$$
such that 
\begin{equation}\label{qII9}
|S^{t_{\nu-1}}(x)| \le 2Mq_{\nu-1}^2\,\,\,\,\,\,\,\, \forall\, x \in J_{\nu}.
\end{equation}
The main difference is that now the right hand side of (\ref{qII9}) does not depend on $t_{\nu-1}$. 
The procedure is the same: for given
  $q_1\le t_1<...<q_{\nu-1}<t_{\nu-1}< q_\nu$ 
we prove the existence of $ t_n, q_{\nu+1}$ and hence $J_{\nu+1}$
taking into account bound (\ref{qII7}) and continuity.
 At the end of the proof instead of  (\ref{tn2}) and (\ref{dzs})  we take $t_{\nu}$ large enough to satisfy
$$
\max_{x\in J_{\nu+1}} |S^{t_\nu}(x)| \le 2Mq_{\nu}^2 
= o (\psi(t_\nu)).
$$
The proof is completed.$\Box$

   \vskip+0.3cm
   
    {\bf 8. About smooth functions.}
    
        \vskip+0.3cm

    In this section,  we discuss examples of smooth functions $f$ for which  
    $
S_{{f}}^{t}(\pmb{\alpha})\to \infty
$
can take place
    (Subsection 8.1) and formulate conditions on smoothness enough to ensure  recurrence  of the sum (\ref{summ0})
    uniformly in $\pmb{x}$  (Section 8.2).
    The conditions on smoothness of function $f$ we formulate in terms of Diophantine exponents 
    ${\omega}(\pmb{\alpha}) ,\hat{\omega}(\pmb{\alpha}),  \hat{\lambda}(\pmb{\alpha}) $ introduced in Section 4.

    \vskip+0.3cm

 {\bf 8.1. Examples of smooth $f$ with $S_f^{t}(\pmb{\alpha})\to \infty$.}
 
     \vskip+0.3cm

We 
follow Kochergin's argument from the beginning of the proof of Theorem 1. But now  instead of functions $F_\nu$ defined in   (\ref{funkktion}) we consider
functions
$$
\underline{F}_\nu (\pmb{x}) =   
\frac{(\sin\, \pi (\pmb{m}_\nu.\pmb{x}))^2 }{ \zeta_\nu} \ge 0,
$$
and
$$
\underline{f}_\nu(\pmb{x}) = \underline{F}_\nu(\pmb{x}+\pmb{\alpha}) - \underline{F}_\nu(\pmb{x})
=
 \frac{ \sin \pi \zeta_\nu}{ \zeta_\nu} \cdot \sin \pi ( 2\pmb{m}_\nu. \pmb{x} +\zeta_\nu),\,\,\,\,\,
\sup_{\nu\in \mathbb{Z}_+}\,\max_{\pmb{x}\in \mathbb{R}^n} |\underline{f}_\nu(\pmb{x}) |\le \pi.
$$
For  real $d\ge 1$ define
$$
\underline{F}(\pmb{x}) = \sum_{\nu=1}^\infty \frac{\underline{F}_\nu (\pmb{x})}{M_\nu^d} 
$$
and
\begin{equation}\label{smi1}
    \underline{f}(\pmb{x}) = \sum_{\nu=1}^\infty \frac{\underline{f}_\nu(\pmb{x})}{M_\nu^d} = 
    \sum_{\nu=1}^\infty  \frac{1}{M_\nu^d} \cdot\frac{ \sin \pi \zeta_\nu}{ \zeta_\nu} \cdot \sin \pi(2\pmb{m}_\nu. \pmb{x} +\zeta_\nu)
\end{equation}
It is clear from the definition (\ref{smi1}) and (\ref{serri}) that function  $    \underline{f}(\pmb{x}) $ has zero mean value and  is $([d]-1)$-smooth.
Analogously to (\ref{sss}) we have  equality
\begin{equation}\label{sss1}
S_{\underline{f}}^{t}(\pmb{\alpha})=
\sum_{k=0}^{t-1} \underline{f}(k\pmb{\alpha }) = \underline{F}(t\pmb{\alpha }). 
\end{equation}
 
   \vskip+0.3cm

{\bf 
Theorem 4.} {\it  Let Diophantine exponents for completely irrational $\pmb{\alpha}$ satisfy
\begin{equation}\label{barbi}
d< 2\hat{\omega} (\pmb{\alpha})- \omega (\pmb{\alpha}).
\end{equation}
Then for the function $\underline{f}$ defined in (\ref{smi1}) one has
$
S_{\underline{f}}^{t}(\pmb{\alpha})\to \infty
$
when $t \to \infty$.
}

   \vskip+0.3cm

   Proof.
   For $ \frac{1}{2\zeta_{\nu-1}}\le t< \frac{1}{2\zeta_{\nu}}$   by (\ref{sss1}) and the  properties (\ref{be1}) of exponents $\omega, \hat{\omega}$ we have
\begin{equation}\label{proof}
S_{\underline{f}}^{t}(\pmb{\alpha})\ge \frac{\underline{F}_\nu(t\pmb{\alpha})}{M_\nu^d} \ge
\frac{\zeta_\nu \cdot (2 t)^2}{M_\nu^d} \ge
\frac{\zeta_\nu }{\zeta_{\nu-1}^2M_\nu^d} \ge  M_\nu^{2\hat{\omega}- \omega- d - \delta}
\end{equation}
for any given $\delta >0$ for $\nu $ large enough. By (\ref{barbi})  the right hand side here tends to infinity, for small $\delta$.$\Box$
  
    \vskip+0.3cm
   {\bf 
Remark 1.} {\it 
If $ \omega = \hat{\omega} = n $ condition (\ref{barbi}) turns into $ d< n$. In fact, if $\pmb{\alpha}$ is badly approximable, even
under assumption 
$d=n$, due to inequalities (\ref{badd}), (\ref{serri}) and (\ref{proof}) we get
$$
\liminf_{t\to \infty}
S_{\underline{f}}^{t}(\pmb{\alpha})>0,
$$
meanwhile function $ \underline{f}$ is $(d-1)$-smooth. This is just the example from \cite{Mnon}.
} 
   
      \vskip+0.3cm
   {\bf  Remark 2.} {\it  If we take completely irrational $\pmb{\alpha}\in \mathbb{R}^n$ with $2\hat{\omega} (\pmb{\alpha})- \omega (\pmb{\alpha})=\frak{g}_n^*$
   where $ \frak{g}_n^*$ is defined in Proposition 1, function $\underline{f}$  from Theorem 4 will be $(w-1)$-smooth for any  integer
   $ w < \frak{g}_n^*\sim \frac{2^{n+1}}{en}, \, n\to \infty$. This gives an improvement of a result from \cite{kon} mentioned in Section 2, {\bf F}.}

    \vskip+0.3cm

 {\bf 8.2 Uniform bounds for $S^t_f (\pmb{\alpha},\pmb{x})$.}
 
     \vskip+0.3cm 
     We mentioned in Section 3 {\bf C} that Sidorov's theorem for absolutely continuous functions  cannot be generalised to  the case $n\ge2$. Indeed   it is shown  in \cite{MK} that  for any  function
    \begin{equation}\label{fufu}
f(\pmb{x}) =  \sum_{\pmb{m}\in \mathbb{Z}^n\setminus\{\pmb{0}\}}
\,\,
f_{\pmb{m}} \exp (2\pi i  \pmb{m}.\pmb{x}),\,\,\,\,\,\,\ n\ge 2,
\end{equation}
such that the set of non-zero Fourier coefficients
$$
\frak{S}_f = \{\pmb{m}\in \mathbb{Z}^n : \,\,\,\, f_{\pmb{m}} \neq 0\} \subset  \mathbb{Z}^n
$$
is not contained in a union of a ball of finite radius and one-dimensional linear  subspace, for any sequence $\sigma_t \downarrow 0$, there exists completely irrational $\pmb{\alpha}\in \mathbb{R}^n$ such that 
 $$
 \sup_{x\in \mathbb{T}^n} \left|
S^t_f (\pmb{\alpha},\pmb{x})\right|  \ge t \cdot \sigma_t
 $$
 for all integer $t $ large enough.

     So  for $n \ge 2$ it is not possible to prove a statement of the form
     \begin{equation}\label{from}
\liminf_{t\to \infty} \sup_{x\in \mathbb{T}^n} \left|
S^t_f (\pmb{\alpha},\pmb{x})\right| = 0.
\end{equation}
for function $f$ of fixed smoothness without additional assumptions on Diophantine properties of $\pmb{\alpha}$. Now we formulate a statement which ensures (\ref{from}) under certain Diophantine conditions on $\pmb{\alpha}$ and conditions  on smoothness of $f$.
     In fact, this is a new version of Proposition 4 from \cite{mou} which involves Diophantine exponents  $\hat{\omega}(\pmb{\alpha})$ and $  \hat{\lambda}(\pmb{\alpha}) $.

   \vskip+0.3cm
   
   {\bf 
Theorem 5.} {\it 
Suppose that 
\begin{equation}\label{rtrt}
\gamma > n+ \frac{\hat{\omega}(\pmb{\alpha})}{\hat{\lambda}(\pmb{\alpha})}
\end{equation}
and in Fourier series  (\ref{fufu})
coefficients are bounded by
\begin{equation}\label{ffff}
 |f_{\pmb{m}} | \le C \left(\max_{1\le j \le n} |m_j|\right)^{-\gamma},\,\,\,\, C, \gamma >0
\end{equation}
(in particular, this means that $f$ is $w$-smooth for any $w<\gamma-n$).
Then
 (\ref{from}) holds.

}

   \vskip+0.3cm
  In particular, (\ref{rtrt}) will be satisfied if we assume
  $$
  \gamma >n+ (n-1)\frac{\hat{\omega}(\pmb{\alpha})^2}{\hat{\omega}(\pmb{\alpha})-1}.
  $$
  This follows from transference inequality
 (\ref{transfer}).
   
      \vskip+0.3cm
  
   {\bf  Remark 3.} {\it In this paper we were interested to formulate the results in terms of classical Diophantine exponents  ${\omega}(\pmb{\alpha})$, $\hat{\omega}(\pmb{\alpha})$ and $  \hat{\lambda}(\pmb{\alpha}) $ and under the condition on the decay of Fourier coefficients of type (\ref{ffff}) which ensures $w$-smoothness.  Traditionally in uniform distribution theory, bounds in terms of product-type approximations are more natural. For example,  assume instead of (\ref{ffff}) another condition
   $$
 |f_{\pmb{m}} | \le C (\overline{\pmb{m}})^{-\gamma_{\text{mult}}},\,\,\,\, C>0,\,\,\, \gamma_{\text{mult}}>0,
 $$
 where 
 $$
 \overline{\pmb{m}}= \prod_{j=1}^n  \max (1,|m_j|)
 $$ 
  and  define multiplicative Diophantine  exponent 
  $$
  \omega_{\text{mult}}(\pmb{\alpha} ) = \sup \{ \gamma \in \mathbb{R}_+:\,\,\,
  \liminf_{\max_j|m_j|\to \infty} (\overline{\pmb{m}})^\gamma \cdot  ||\pmb{m}.\pmb{\alpha}|| <\infty
  \}.
  $$
  Then it is well known that condition $ \gamma_{\text{mult}}>  \omega_{\text{mult}}$ leads to 
 \begin{equation}\label{known}
\sup_{t \in \mathbb{Z}_+} \sup_{x\in \mathbb{T}^n} \left|
S^t_f (\pmb{\alpha},\pmb{x})\right| <\infty,
\end{equation}
and as for almost all $\pmb{\alpha}$ we have $\omega_{\text{mult}}(\pmb{\alpha} ) =1$, we see that  
for almost all $\pmb{\alpha}$ we have (\ref{known}) under the assumption $ \gamma_{\text{mult}}>1$.
 
   }

      \vskip+0.3cm
   Sketch of the proof of Theorem 5.
       The proof is completely similar to that from \cite{mou}.
    The difference is that instead of pigeon hole principle one should take into account that for any $\delta>0 $ and for any $Q$ large enough there exists
    a positive integer $t$ such that
    $$
    1\le t \le Q, \,\,\,\, \max_{1\le j \le n} ||t\alpha_j|| \le Q^{-\hat{\lambda}+\delta}.
    $$
    From Fourier expansion  for such $t$ we get inequalities
    $$
    \sup_{x\in \mathbb{T}^n} \left|
S^t_f (\pmb{\alpha},\pmb{x})\right|
\ll
$$
$$
\ll
\frac{Q^{-\hat{\lambda}+\delta}}{\zeta_\nu}
\sum_{\begin{array}{c}\pmb{m}\in \mathbb{Z}^n: \cr \max_j |m_j|< M_{\nu+1}\end{array}}
\frac{1}{\left(\max_{j }|m_j|\right)^\gamma
} +
Q\sum_{\begin{array}{c}  \pmb{m}\in \mathbb{Z}^n: \cr \max_j |m_j|\ge M_{\nu+1} \end{array}}
 \frac{1}{\left(\max_{j }|m_j|\right)^\gamma}
 \ll
 $$
 $$
 \ll
 Q^{-\hat{\lambda}+\delta}\zeta_\nu^{-1} + Q M_{\nu+1}^{n-\gamma},
 $$
 for a proper value of $\nu$ with $ \zeta_\nu \ge M_{\nu+1}^{-\hat{\omega}+\delta}$. Then one should optimise this inequality  by choosing  
$Q = M^{\frac{\gamma- n+ \hat{\omega}}{1+\hat{\lambda}} + \delta'}$.
 We leave the details to the reader . $\Box$

    \vskip+0.3cm
{\bf Acknowledgements.}
 
      \vskip+0.3cm
The author is extremely grateful to Anton Shutov for many  valuable comments on manuscript and 
important suggestions. The author also is very grateful  to the anonymous referee for very careful reading of the manuscript and valuable comments.
  \vskip+0.3cm
The author is a winner of the “Leader Researcher” contest conducted by Theoretical Physics and Mathematics Advancement Foundation “BASIS” and would like to thank its sponsors and jury.

The author thanks 
Yau Mathematical Sciences Center (YMSC) 
for hospitality 
and wonderful opportunity to work in Tsinghua University.

  \vskip+0.3cm
This research is supported by joint FWF-Projekt I 5554 and RSF project 22-41-05001 \\ https://rscf.ru/en/project/22-41-05001/).

      \vskip+0.3cm

 \vskip+1cm
 
 Nikolay Moshchevitin,
 
 Technische Universität Wien
 
 and
 
 Faculty of Computer Science, HSE University,
 
Moscow Center of Fundamental and Applied Mathematics

 \vskip+0.3cm
    
nikolai.moshchevitin@tuwien.ac.at,  nikolaus.moshchevitin@gmail.com
 \end{document}